\documentclass{amsart}
\title{The Carmichael numbers up to $10^{16}$}
\author{Richard G.E. Pinch}
\address{Queens' College, Silver Street, Cambridge CB3 9ET, U.K.}
\email{rgep@cam.ac.uk}

\date{18 March 1998}

\usepackage{amsmath}
\usepackage{amssymb}
\usepackage[dvips]{epsfig}

\newcommand{\divides}{\vert}
\newcommand{\congruent}{\equiv}

\def\abs|#1|{{\left\vert{#1}\right\vert}}
\def\O(#1){\mbox{O}{\left({#1}\right)}}
\def\C(#1){C{\left({#1}\right)}}
\def\GF(#1){{\mathbb{F}_{#1}}}			
\def\gen<#1>{{\left\langle{#1}\right\rangle}}
\def\paren(#1){{\left({#1}\right)}}

\begin{document}

\begin{abstract}
We extend our previous computations to show that there are
246683 Carmichael numbers up to $10^{16}$.  As before,
the numbers were generated by a back-tracking search for
possible prime factorisations together with a
``large prime variation''.  We present further statistics on
the distribution of Carmichael numbers.
\end{abstract}

\maketitle

\suppressfloats

\section{Introduction}

A {\em Carmichael number} $N$ is a composite number $N$ with the property
that for every $b$ prime to $N$ we have $b^{N-1} \congruent 1 \mod N$.
It follows that a Carmichael number $N$
must be square-free, with at least three prime factors, and
that $p-1 \divides N-1$ for every prime $p$ dividing $N$:
conversely, any such $N$ must be a Carmichael number.

For background on Carmichael numbers and details of previous computations
we refer to our previous paper \cite{Pin:car15}: in that paper we described
the computation of the Carmichael numbers up to $10^{15}$ and presented
some statistics.  These computations have since been extended to $10^{16}$,
using the same techniques, and we present further statistics.

The complete list of Carmichael numbers up to $10^{16}$
is available by anonymous FTP
from {\tt ftp.dpmms.cam.ac.uk} in directory {\tt /pub/Carmichael}.

\section{Statistics}

\begin{table}[phtb]
\begin{tabular}{||r|r||}
\hline
 $n$  & $C\left(10^n\right)$  \\
\hline
 3  &      1  \\
 4  &      7  \\
 5  &     16  \\
 6  &     43  \\
 7  &    105  \\
 8  &    255  \\
 9  &    646  \\
 10  &   1547  \\
 11  &   3605  \\
 12  &   8241  \\
 13  &  19279  \\
 14  &  44706  \\
 15  & 105212 \\
 16  & 246683 \\
\hline
\end{tabular}
\caption{Distribution of Carmichael numbers up to $10^{16}$.}
\label{table0}
\end{table}

\begin{table}[phtb]
\begin{tabular}{||r|r|r|r|r|r|r|r|r|r||}
\hline
 $X$  &  3  &  4  &  5  &  6  &  7  &  8  &  9  &  10  &  total  \\
\hline
   3  &     1  &     0  &     0  &     0  &     0  &     0  &     0  &  0  &      1  \\
   4  &     7  &     0  &     0  &     0  &     0  &     0  &     0  &  0  &      7  \\
   5  &    12  &     4  &     0  &     0  &     0  &     0  &     0  &  0  &     16  \\
   6  &    23  &    19  &     1  &     0  &     0  &     0  &     0  &  0  &     43  \\
   7  &    47  &    55  &     3  &     0  &     0  &     0  &     0  &  0  &    105  \\
   8  &    84  &   144  &    27  &     0  &     0  &     0  &     0  &  0  &    255  \\
   9  &   172  &   314  &   146  &    14  &     0  &     0  &     0  &  0  &    646  \\
  10  &   335  &   619  &   492  &    99  &     2  &     0  &     0  &  0  &   1547  \\
  11  &   590  &  1179  &  1336  &   459  &    41  &     0  &     0  &  0  &   3605  \\
  12  &  1000  &  2102  &  3156  &  1714  &   262  &     7  &     0  &  0  &   8241  \\
  13  &  1858  &  3639  &  7082  &  5270  &  1340  &    89  &     1  &  0  &  19279  \\
  14  &  3284  &  6042  & 14938  & 14401  &  5359  &   655  &    27  &  0  &  44706  \\
  15  &  6083  &  9938  & 29282  & 36907  & 19210  &  3622  &   170  &  0  & 105212  \\
  16  & 10816  & 16202  & 55012  & 86696  & 60150  & 16348  &  1436  & 23  & 246683  \\
\hline
\end{tabular}
\caption{Values of $\C(X)$ and $\C(d,X)$ for $d \le 10$ and $X$ in
powers of 10 up to $10^{16}$.}
\label{table1}
\end{table}

We have shown that there are 246683 Carmichael numbers up to $10^{16}$,
all with at most 10 prime factors.  We let $\C(X)$ denote
the number of Carmichael numbers less than $X$ and $\C(d,X)$ denote the
number with exactly $d$ prime factors.  Table \ref{table0} gives the
values of $\C(X)$ and Table \ref{table1} the values
of $\C(d,X)$ for $X$ in powers of 10 up to $10^{16}$.

We have used the same methods to calculate the smallest Carmichael numbers
with $d$ prime factors for $d$ up to 20.  The results are given in
Table \ref{table2}.

\begin{table}[phtb]
\begin{tabular}{||r|l||}
\hline
  $d$  &  $N$  \\
  {}  &  \mbox{factors}  \\
\hline
   3  &  \hfill                561  \\
  {}  &   3.11.17  \\
\hline
   4  &  \hfill               41041  \\
  {}  &   7.11.13.41  \\
\hline
   5  &  \hfill              825265  \\
  {}  &   5. 7.17.19.73  \\
\hline
   6  &  \hfill           321197185  \\
  {}  &   5.19.23.29.37.137  \\
\hline
   7  &  \hfill          5394826801  \\
  {}  &   7.13.17.23.31.67.73  \\
\hline
   8  &  \hfill        232250619601  \\
  {}  &   7.11.13.17.31.37.73.163  \\
\hline
   9  &  \hfill       9746347772161  \\
  {}  &   7.11.13.17.19.31.37.41.641  \\
\hline
  10  &  \hfill    1436697831295441  \\
  {}  &  11.13.19.29.31.37.41.43.71.127  \\
\hline
  11  &  \hfill   60977817398996785  \\
  {}  &   5. 7.17.19.23.37.53.73.79.89.233  \\
\hline
  12  &  \hfill 7156857700403137441  \\
  {}  &  11.13.17.19.29.37.41.43.61.97.109.127  \\
\hline
  13  &  \hfill   1791562810662585767521  \\
  {}  &  11.13.17.19.31.37.43.71.73.97.109.113.127  \\
\hline
  14  &  \hfill  87674969936234821377601  \\
  {}  &   7.13.17.19.23.31.37.41.61.67.89.163.193.241  \\
\hline
  15  &  \hfill 6553130926752006031481761  \\
  {}  &  11.13.17.19.29.31.41.43.61.71.73.109.113.127.181  \\
\hline
  16  &  \hfill 1590231231043178376951698401  \\
  {}  &  17.19.23.29.31.37.41.43.61.67.71.73.79. 97.113.199  \\
\hline
  17  &  \hfill 35237869211718889547310642241  \\
  {}  &  13.17.19.23.29.31.37.41.43.61.67.71.73. 97.113.127.211  \\
\hline
  18  &  \hfill 32809426840359564991177172754241  \\
  {}  &  13.17.19.23.29.31.37.41.43.61.67.71.73. 97.127.199.281.397  \\
\hline
  19  &  \hfill 2810864562635368426005268142616001  \\
  {}  & 13.17.19.23.29.31.37.41.43.61.67.71.73.109.113.127.151.281.353  \\
\hline
  20  &  \hfill 349407515342287435050603204719587201  \\
  {}  &  11.13.17.19.29.31.37.41.43.61.71.73.97.101.109.113.151.181.193.641 \\
\hline
\end{tabular}
\caption{The smallest Carmichael numbers with $d$ prime factors for $d$ up to 20.}
\label{table2}
\end{table}

\begin{figure}
\epsfig{file=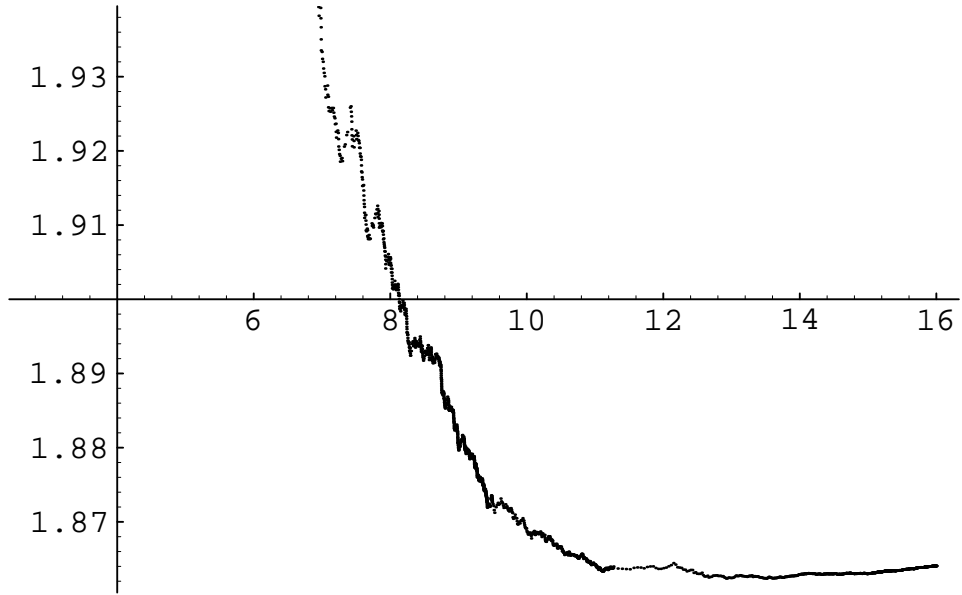,width=\textwidth}
\vspace{-1truecm}
\caption{$k(X)$ versus $X$ (expressed in powers of 10).}
\label{figure1}
\end{figure}

\begin{table}[phtb]
\begin{tabular}{||r|r|c||}
\hline
   $n$  &  $k\left(10^n\right)$  &  $\C(10^n)/\C(10^{n-1})$ \\
\hline
    3  &  2.93319  &        \\
    4  &  2.19547  &  7.000 \\
    5  &  2.07632  &  2.286 \\
    6  &  1.97946  &  2.688 \\
    7  &  1.93388  &  2.441 \\
    8  &  1.90495  &  2.429 \\
    9  &  1.87989  &  2.533 \\
   10  &  1.86870  &  2.396 \\
   11  &  1.86421  &  2.330 \\
   12  &  1.86377  &  2.286 \\
   13  &  1.86240  &  2.339 \\
   14  &  1.86293  &  2.319 \\
   15  &  1.86301  &  2.353 \\
   16  &  1.86406  &  2.335 \\
\hline
\end{tabular}
\caption{The function $k(X)$ and growth of $\C(X)$ for $X = 10^n$,
$n \le 16$.}
\label{table3}
\end{table}

In Table \ref{table3} and Figure \ref{figure1} we tabulate the function
$k(X)$, defined by Pomerance, Selfridge and Wagstaff \cite{PSW:pseudo} by
$$
\C(X) = X \exp\left(-k(X) \frac{\log X \log\log\log X}{\log\log X}\right) .
$$
They proved that $\liminf k \ge 1$
and suggested that $\limsup k$ might be $2$, although they also observed
that within the range of their tables $k(X)$ is decreasing:
Pomerance \cite{Pom:distpsp},\cite{Pom:2meth}
gave a heuristic argument suggesting that $\lim k = 1$.
The decrease in $k$ is reversed between $10^{13}$ and $10^{14}$: see
Figure \ref{figure1}.  We find no clear support from our computations
for any conjecture on a limiting value of $k$.

In Table \ref{table3} we also give the ratios $\C(10^n) / \C(10^{n-1})$
investigated by Swift \cite{Swi:car}.
Swift's ratio, again initially decreasing, also increases again before
$10^{15}$.

\begin{figure}
\epsfig{file=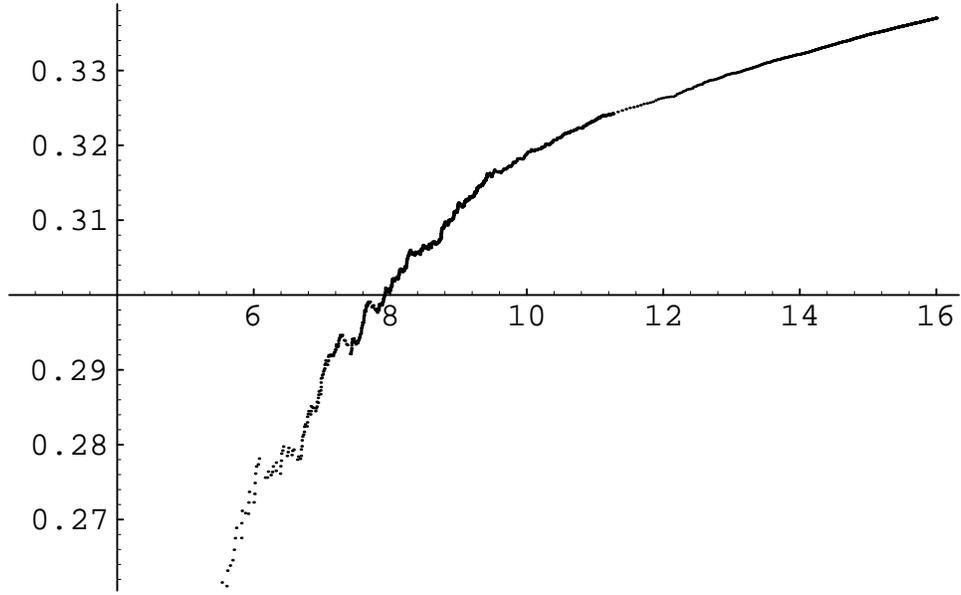,width=\textwidth}
\vspace{-1truecm}
\caption{$C(X)$ as a power of $X$ (expressed in powers of 10).}
\label{figure2}
\end{figure}

\begin{table}[phtb]
\begin{tabular}{||l|c||}
\hline
 $n$  &  $\log C\left(10^n\right) / (n \log 10)$ \\
\hline
 4   &  .21127  \\
 5   &  .24082  \\
 6   &  .27224  \\
 7   &  .28874  \\
 8   &  .30082  \\
 9   &  .31225  \\
 {10}  &  .31895  \\
 {11}  &  .32336  \\
 {12}  &  .32633  \\
 {13}  &  .32962  \\
 {14}  &  .33217  \\
 {15}  &  .33480  \\
 {16}  &  .33700  \\
\hline
\end{tabular}
\caption{$\C(X)$ as a power of $X$.}
\label{table00}
\end{table}

In Table \ref{table00} and Figure \ref{figure2} we see that within the
range of our computations $C(X)$ is a slowly growing power of $X$:
about $X^{0.337}$ for $X$ around $10^{16}$.

In Table \ref{table4} we give the number of Carmichael numbers in each
class modulo $m$ for $m = 5$, 7, 11 and 12.

In Tables \ref{table5} and \ref{table6} we give the number of Carmichael
numbers divisible by primes $p$ up to 97.  In Table \ref{table5} we count
all Carmichael numbers divisible by $p$: in Table \ref{table6} we count
only those for which $p$ is the smallest prime factor.

The largest prime factor of a Carmichael number up
to $10^{16}$ is 68786257, dividing
$$
9463098235353841 = 13 \cdot 31 \cdot 541 \cdot 631 \cdot 68786257
$$
and the largest prime to occur as the smallest prime factor of a Carmichael
number in this range is 174763, dividing
$$
9585921133193329 = 174763 \cdot 199729 \cdot 274627 .
$$
We note that this number is of the form $(7k+1)(8k+1)(11k+1)$
with $k = 24966$.

\begin{table}[phtb]
\begin{tabular}{||r|r|r|r|r|r|r|r|r||}
\hline
 $m$ & $c$ & $25.10^9$ & $10^{11}$ & $10^{12}$ & $10^{13}$ & $10^{14}$ & $10^{15}$ & $10^{16}$ \\
\hline
    5   &  0  &    203  &  312  &  627  &  1330  &   2773  &   5814  & 12200	\\
        &  1  &   1652  & 2785  & 6575  & 15755  &  37467  &  90167  & 215713	\\
        &  2  &     82  &  154  &  327  &   702  &   1484  &   3048  & 6094	\\
        &  3  &    102  &  172  &  344  &   725  &   1463  &   3059  & 6285	\\
        &  4  &    124  &  182  &  368  &   767  &   1519  &   3124  & 6391	\\
\hline
    7   &  0  &    401  &  634  & 1334  &   2774  &   5891  &  12691  & 27550	\\
        &  1  &   1096  & 1885  & 4613  &  11447  &  28001  &  69131  & 168856	\\
        &  2  &    105  &  186  &  432  &    967  &   2109  &   4599  & 10011	\\
        &  3  &    152  &  232  &  496  &   1055  &   2178  &   4707  & 10039	\\
        &  4  &    129  &  211  &  450  &    985  &   2122  &   4592  & 9944	\\
        &  5  &    138  &  222  &  454  &   1033  &   2224  &   4777  & 10125	\\
        &  6  &    142  &  235  &  462  &   1018  &   2181  &   4715  & 10158	\\
\hline
   11   &  0  &    335  &  547  &  1324  &  3006  &   7032  &  16563  & 38576	\\
        &  1  &    640  & 1131  &  2770  &  6786  &  16548  &  40891  & 100071	\\
        &  2  &    139  &  217  &  473  &   1068  &   2361  &   5338  & 12054	\\
        &  3  &    142  &  220  &  457  &   1045  &   2348  &   5319  & 12186	\\
        &  4  &    104  &  187  &  442  &   1026  &   2317  &   5261  & 11917	\\
        &  5  &    152  &  243  &  466  &   1066  &   2370  &   5316  & 12194	\\
        &  6  &    116  &  198  &  440  &   1061  &   2400  &   5384  & 12155	\\
        &  7  &    122  &  195  &  458  &   1023  &   2223  &   5165  & 11853	\\
        &  8  &    129  &  222  &  475  &   1107  &   2450  &   5449  & 12012	\\
        &  9  &    131  &  218  &  465  &   1042  &   2285  &   5179  & 11835	\\
       &  10  &    153  &  227  &  471  &   1049  &   2372  &   5347  & 11830	\\
\hline
   12   &  1  &   2071  & 3462  & 7969  & 18761  & 43760  &  103428  & 243382	\\
        &  3  &      0  &    0  &    1  &     2  &     2  &       5  & 5	\\
        &  5  &     20  &   32  &   64  &   124  &   228  &     448  & 805	\\
        &  7  &     47  &   75  &  147  &   289  &   547  &    1027  & 1906	\\
        &  9  &     25  &   36  &   60  &   103  &   165  &     294  & 560	\\
        & 11  &      0  &    0  &    0  &     0  &     4  &      10  & 25	\\
\hline
\end{tabular}
\caption{The number of Carmichael numbers in each class modulo $m$
for $m = 5$, 7, 11 and 12.}
\label{table4}
\end{table}

\begin{table}[phtb]
\begin{tabular}{||r|r|r|r|r|r|r|r||}
\hline
$p$ & $25.10^9$ & $10^{11}$ & $10^{12}$ & $10^{13}$ & $10^{14}$ & $10^{15}$ & $10^{16}$ \\
\hline
     3  &   25  &    36  &    61  &   105  &   167  &   299  &    565 \\
     5  &  203  &   312  &   627  &  1330  &  2773  &  5814  &  12200 \\
     7  &  401  &   634  &  1334  &  2774  &  5891  & 12691  &  27550 \\
    11  &  335  &   547  &  1324  &  3006  &  7032  & 16563  &  38576 \\
    13  &  483  &   807  &  1784  &  3998  &  9045  & 20758  &  47785 \\
    17  &  293  &   489  &  1182  &  2817  &  6640  & 16019  &  38302 \\
    19  &  372  &   608  &  1355  &  3345  &  7797  & 18638  &  44389 \\
    23  &  113  &   207  &   507  &  1282  &  3135  &  7716  &  18867 \\
    29  &  194  &   336  &   832  &  2094  &  5158  & 12721  &  31110 \\
    31  &  335  &   571  &  1320  &  3086  &  7270  & 17382  &  41440 \\
    37  &  320  &   535  &  1270  &  2926  &  6826  & 16220  &  38647 \\
    41  &  227  &   390  &  1001  &  2418  &  5896  & 14344  &  34759 \\
    43  &  184  &   296  &   772  &  1920  &  4663  & 11594  &  28650 \\
    47  &   53  &    80  &   199  &   492  &  1223  &  2873  &   6810 \\
    53  &   92  &   160  &   351  &   813  &  2041  &  5143  &  12256 \\
    59  &   26  &    41  &    92  &   262  &   644  &  1611  &   3959 \\
    61  &  269  &   453  &  1075  &  2542  &  6047  & 14429  &  34503 \\
    67  &  110  &   178  &   407  &  1063  &  2540  &  6306  &  15295 \\
    71  &  104  &   194  &   521  &  1320  &  3351  &  8546  &  21485 \\
    73  &  198  &   348  &   849  &  2145  &  4925  & 11929  &  29072 \\
    79  &   64  &   107  &   247  &   686  &  1728  &  4318  &  10693 \\
    83  &   14  &    24  &    56  &   137  &   340  &   838  &   1929 \\
    89  &   68  &   131  &   320  &   788  &  1951  &  4981  &  12178 \\
    97  &  123  &   193  &   495  &  1277  &  3123  &  7594  &  18706 \\
\hline
\end{tabular}
\caption{Primes occurring in Carmichael numbers.}
\label{table5}
\end{table}

\begin{table}[phtb]
\begin{tabular}{||r|r|r|r|r|r|r|r||}
\hline
$p$ & $25.10^9$ & $10^{11}$ & $10^{12}$ & $10^{13}$ & $10^{14}$ & $10^{15}$ & $10^{16}$ \\	
\hline
     3  &   25  &    36  &    61  &   105  &   167  &   299  &    565  \\
     5  &  202  &   309  &   624  &  1325  &  2765  &  5797  &  12175  \\
     7  &  364  &   579  &  1218  &  2557  &  5461  & 11874  &  25915  \\
    11  &  263  &   428  &  1071  &  2509  &  5979  & 14397  &  33893  \\
    13  &  237  &   431  &  1058  &  2462  &  5699  & 13514  &  32025  \\
    17  &  117  &   206  &   496  &  1318  &  3244  &  8114  &  20206  \\
    19  &  152  &   244  &   532  &  1401  &  3358  &  8141  &  20020  \\
    23  &   37  &    78  &   207  &   535  &  1360  &  3317  &   8195  \\
    29  &   55  &   103  &   284  &   729  &  1822  &  4659  &  11577  \\
    31  &  101  &   168  &   390  &   876  &  2116  &  5153  &  12575  \\
    37  &   60  &    95  &   219  &   551  &  1401  &  3418  &   8594  \\
    41  &   35  &    68  &   171  &   414  &  1092  &  2736  &   6788  \\
    43  &   35  &    65  &   168  &   403  &   943  &  2308  &   5520  \\
    47  &   14  &    16  &    36  &    81  &   195  &   459  &   1135  \\
    53  &   19  &    30  &    55  &   147  &   363  &   973  &   2327  \\
    59  &    2  &     4  &    11  &    43  &   100  &   272  &    618  \\
    61  &   34  &    58  &   148  &   364  &   851  &  1978  &   4722  \\
    67  &    8  &    18  &    50  &   123  &   317  &   815  &   1950  \\
    71  &   15  &    25  &    66  &   161  &   389  &   979  &   2480  \\
    73  &   14  &    28  &    68  &   175  &   406  &  1015  &   2508  \\
    79  &    4  &    10  &    17  &    66  &   175  &   467  &   1163  \\
    83  &    1  &     1  &     4  &     8  &    39  &    79  &    175  \\
    89  &   10  &    16  &    23  &    55  &   148  &   409  &   1003  \\
    97  &   10  &    20  &    50  &   106  &   261  &   606  &   1413  \\
\hline
\end{tabular}
\caption{Primes occurring as least prime factor in Carmichael numbers.}
\label{table6}
\end{table}
	
\clearpage


\ifx\undefined\bysame
\newcommand{\bysame}{\leavevmode\hbox to3em{\hrulefill}\,}
\fi


\begin{thebibliography}{1}

\bibitem{Mol:ntapps}
R.A. Mollin (ed.), {\em Number theory and its applications}, Dordrecht, Kluwer
  Academic, 1989, Proceedings of the NATO Advanced Study Institute on Number
  Theory and Applications.

\bibitem{Pin:car15}
R.G.E. Pinch, {\em The {Carmichael} numbers up to $10^{15}$}, Math. Comp. {\bf
  61} (1993), 381--391, {Lehmer} memorial issue.

\bibitem{Pom:distpsp}
C.~Pomerance, {\em On the distribution of pseudoprimes}, Math. Comp. {\bf 37}
  (1981), 587--593.

\bibitem{Pom:2meth}
\bysame, {\em Two methods in elementary analytic number theory}, In Mollin
  \cite{Mol:ntapps}, Proceedings of the NATO Advanced Study Institute on Number
  Theory and Applications.

\bibitem{PSW:pseudo}
C.~Pomerance, J.L. Selfridge, and S.S. Wagstaff~jr, {\em The pseudoprimes up to
  $25.10^9$}, Math. Comp. {\bf 35} (1980), no.~151, 1003--1026.

\bibitem{Swi:car}
J.D. Swift, {\em Review 13[9] --- table of {Carmichael} numbers to $10^9$},
  Math. Comp. {\bf 29} (1975), 338--339.

\end{thebibliography}
\end{document}